\documentclass[final]{siamltex}
\pdfoutput=1
\usepackage{amsmath}
\usepackage{amssymb}
\usepackage{epsfig}
\usepackage{amsbsy}

\newcommand{\rk}{\mathrm{rank}}
\newcommand{\s}{\mathrm{support}}
\newcommand{\As}{A^\star}
\newcommand{\Bs}{B^\star}
\newcommand{\R}{\mathbb{R}}

\newcommand{\inc}{\mathrm{inc}}
\newcommand{\degr}{\mathrm{deg}}

\newtheorem{THEO}{Theorem}
\newtheorem{CORL}[THEO]{Corollary}
\newtheorem{LEMM}{Lemma}
\newtheorem{PROP}{Proposition}

\title{Rank-Sparsity Incoherence for Matrix Decomposition}

\author{Venkat Chandrasekaran, Sujay Sanghavi, Pablo A.  Parrilo,\and
Alan S. Willsky\thanks{Venkat Chandrasekaran, Pablo A.  Parrilo, and
Alan S. Willsky are with the Laboratory for Information and Decision
Systems, Department of Electrical Engineering and Computer Science,
Massachusetts Institute of Technology, Cambridge, MA 02139 ({\tt
venkatc@mit.edu; parrilo@mit.edu; willsky@mit.edu}).  Sujay Sanghavi
is with the School of Electrical and Computer Engineering, Purdue
University, West Lafayette, IN 47907 ({\tt
sanghavi@ecn.purdue.edu}). This work was supported by MURI AFOSR
grant FA9550-06-1-0324, MURI AFOSR grant FA9550-06-1-0303, and NSF
FRG 0757207. A preliminary version of this work will appear at the
15th IFAC Symposium on System Identification \cite{Cha:09}.
Submitted on June 11, 2009.}}

\begin{document}

\maketitle

\begin{abstract}
Suppose we are given a matrix that is formed by adding an unknown
sparse matrix to an unknown low-rank matrix. Our goal is to
decompose the given matrix into its sparse and low-rank components.
Such a problem arises in a number of applications in model and
system identification, and is NP-hard in general. In this paper we
consider a convex optimization formulation to splitting the
specified matrix into its components, by minimizing a linear
combination of the $\ell_1$ norm and the nuclear norm of the
components. We develop a notion of \emph{rank-sparsity incoherence},
expressed as an uncertainty principle between the sparsity pattern
of a matrix and its row and column spaces, and use it to
characterize both fundamental identifiability as well as
(deterministic) sufficient conditions for exact recovery. Our
analysis is geometric in nature, with the tangent spaces to the
algebraic varieties of sparse and low-rank matrices playing a
prominent role. When the sparse and low-rank matrices are drawn from
certain natural random ensembles, we show that the sufficient
conditions for exact recovery are satisfied with high probability.
We conclude with simulation results on synthetic matrix
decomposition problems.

\end{abstract}



\pagestyle{myheadings} \thispagestyle{plain} \markboth{V.
Chandrasekaran, S. Sanghavi, P. A. Parrilo, and A. S.
Willsky}{Rank-Sparsity Incoherence for Matrix Decomposition}

\section{Introduction}
\label{sec:intro}


Complex systems and models arise in a variety of problems in science
and engineering. In many applications such complex systems and
models are often composed of multiple simpler systems and models.
Therefore, in order to better understand the behavior and properties
of a complex system a natural approach is to decompose the system
into its simpler components. In this paper we consider matrix
representations of systems and statistical models in which our
matrices are formed by adding together \emph{sparse} and
\emph{low-rank} matrices. We study the problem of recovering the
sparse and low-rank components given no prior knowledge about the
sparsity pattern of the sparse matrix, or the rank of the low-rank
matrix. We propose a tractable convex program to recover these
components, and provide sufficient conditions under which our
procedure recovers the sparse and low-rank matrices \emph{exactly}.


Such a decomposition problem arises in a number of settings, with
the sparse and low-rank matrices having different interpretations
depending on the application. In a statistical model selection
setting, the sparse matrix can correspond to a Gaussian graphical
model \cite{Lau:96} and the low-rank matrix can summarize the effect
of latent, unobserved variables. Decomposing a given model into
these simpler components is useful for developing efficient
estimation and inference algorithms. In computational complexity,
the notion of \emph{matrix rigidity} \cite{Val:77} captures the
smallest number of entries of a matrix that must be changed in order
to reduce the rank of the matrix below a specified level (the
changes can be of arbitrary magnitude). Bounds on the rigidity of a
matrix have several implications in complexity theory \cite{Lok:95}.
Similarly, in a system identification setting the low-rank matrix
represents a system with a small model order while the sparse matrix
represents a system with a sparse impulse response. Decomposing a
system into such simpler components can be used to provide a
simpler, more efficient description.


\subsection{Our results} Formally the decomposition problem we are
interested can be defined as follows:

\paragraph*{Problem} Given $C = \As + \Bs$ where $\As$ is an unknown
sparse matrix and $\Bs$ is an unknown low-rank matrix, recover $\As$
and $\Bs$ from $C$ using no additional information on the sparsity
pattern and/or the rank of the components.

In the absence of any further assumptions, this decomposition
problem is fundamentally ill-posed. Indeed, there are a number of
scenarios in which a unique splitting of $C$ into ``low-rank'' and
``sparse'' parts may not exist; for example, the low-rank matrix may
itself be very sparse leading to identifiability issues. In order to
characterize when such a decomposition is possible we develop a
notion of \emph{rank-sparsity incoherence}, an uncertainty principle
between the sparsity pattern of a matrix and its row/column spaces.
This condition is based on quantities involving the tangent spaces
to the algebraic variety of sparse matrices and the algebraic
variety of low-rank matrices \cite{Har:95}.

Two natural identifiability problems may arise. The first one occurs
if the low-rank matrix itself is very sparse. In order to avoid such a
problem we impose certain conditions on the row/column spaces of the
low-rank matrix.  Specifically, for a matrix $M$ let $T(M)$ be the
tangent space at $M$ with respect to the variety of all matrices with
rank less than or equal to $\rk(M)$.  Operationally, $T(M)$ is the
span of all matrices with row-space contained in the row-space of $M$
or with column-space contained in the column-space of $M$; see
(\ref{eq:t}) for a formal characterization. Let $\xi(M)$ be defined as
follows:
\begin{equation}
\xi(M) \triangleq \max_{N \in T(M), ~ \|N\| \leq 1} \|N\|_\infty.
\label{eq:xim}
\end{equation}
Here $\|\cdot\|$ is the spectral norm (i.e., the largest singular
value), and $\|\cdot\|_\infty$ denotes the largest entry in
magnitude. Thus $\xi(M)$ being small implies that (appropriately
scaled) elements of the tangent space $T(M)$ are ``diffuse'', i.e.,
these elements are not too sparse; as a result $M$ cannot be very
sparse. As shown in Proposition~\ref{prop:lr} (see
Section~\ref{subsec:specialslr}) a low-rank matrix $M$ with row/column
spaces that are not closely aligned with the coordinate axes has small
$\xi(M)$.

The other identifiability problem may arise if the sparse matrix has
all its support concentrated in one column; the entries in this column
could negate the entries of the corresponding low-rank matrix, thus
leaving the rank and the column space of the low-rank matrix
unchanged. To avoid such a situation, we impose conditions on the
sparsity pattern of the sparse matrix so that its support is not too
concentrated in any row/column. For a matrix $M$ let $\Omega(M)$ be
the tangent space at $M$ with respect to the variety of all matrices
with number of non-zero entries less than or equal to $|\s(M)|$. The
space $\Omega(M)$ is simply the set of all matrices that have support
contained within the support of $M$; see~(\ref{eq:om}). Let
$\mu(M)$ be defined as follows:
\begin{equation}
\mu(M) \triangleq \max_{N \in \Omega(M), ~ \|N\|_\infty \leq 1}
\|N\|. \label{eq:mum}
\end{equation}
The quantity $\mu(M)$ being small for a matrix implies that the
\emph{spectrum} of any element of the tangent space $\Omega(M)$ is
``diffuse'', i.e., the singular values of these elements are not too
large. We show in Proposition~\ref{prop:sp} (see
Section~\ref{subsec:specialslr}) that a sparse matrix $M$ with
``bounded degree'' (a small number of non-zeros per row/column) has
small $\mu(M)$.

For a given matrix $M$, it is impossible for both quantities $\xi(M)$
and $\mu(M)$ to be simultaneously small. Indeed, we prove that for any
matrix $M \neq 0$ we must have that $\xi(M) \mu(M) \geq 1$ (see
Theorem~\ref{theo:unp} in Section~\ref{subsec:unp}). Thus, this
\emph{uncertainty principle} asserts that there is no non-zero matrix
$M$ with all elements in $T(M)$ being diffuse \emph{and} all elements
in $\Omega(M)$ having diffuse spectra. As we describe later, the
quantities $\xi$ and $\mu$ are also used to characterize fundamental
identifiability in the decomposition problem.


In general solving the decomposition problem is NP-hard; hence, we
consider tractable approaches employing recently well-studied convex
relaxations. We formulate a convex optimization problem for
decomposition using a combination of the $\ell_1$ norm and the
nuclear norm. For any matrix $M$ the $\ell_1$ norm is given by
\begin{equation*}
\|M\|_1 = \sum_{i,j} |M_{i,j}|,
\end{equation*}
and the nuclear norm, which is the sum of the singular values, is
given by
\begin{equation*}
\|M\|_{\ast} = \sum_{k} \sigma_k(M),
\end{equation*}
where $\{\sigma_k(M)\}$ are the singular values of $M$.  The $\ell_1$
norm has been used as an effective surrogate for the number of
non-zero entries of a vector, and a number of results provide
conditions under which this heuristic recovers sparse solutions to
ill-posed inverse problems \cite{Don:06}. More recently, the nuclear
norm has been shown to be an effective surrogate for the rank of a
matrix \cite{Faz:02}. This relaxation is a generalization of the
previously studied trace-heuristic that was used to recover low-rank
positive semidefinite matrices \cite{Mes:97}. Indeed, several papers
demonstrate that the nuclear norm heuristic recovers low-rank matrices
in various rank minimization problems \cite{Rec:07,Can:08}.  Based on
these results, we propose the following optimization formulation to
recover $\As$ and $\Bs$ given $C = \As + \Bs$:
\begin{equation}
\begin{aligned}
(\hat{A},\hat{B}) = \arg \min_{A,B} &  ~~ \gamma \|A\|_{1} +
\|B\|_{\ast} \\ \mbox{s.t.} & ~~ A + B = C.
\end{aligned}
\label{eq:sdp}
\end{equation}
Here $\gamma$ is a parameter that provides a trade-off between the
low-rank and sparse components. This optimization problem is convex,
and can in fact be rewritten as a semidefinite program (SDP)
\cite{Van:96} (see Appendix~\ref{sec:sdp}).

We prove that $(\hat{A},\hat{B}) = (\As, \Bs)$ is the unique optimum
of (\ref{eq:sdp}) for a range of $\gamma$ if $\mu(\As) \xi(\Bs) <
\frac{1}{6}$ (see Theorem~\ref{theo:main} in
Section~\ref{subsec:mainres}). Thus, the conditions for \emph{exact}
recovery of the sparse and low-rank components via the convex
program (\ref{eq:sdp}) involve the tangent-space-based quantities
defined in (\ref{eq:xim}) and (\ref{eq:mum}). Essentially these
conditions specify that each element of $\Omega(\As)$ must have a
diffuse spectrum, \emph{and} every element of $T(\Bs)$ must be
diffuse. In a sense that will be made precise later, the condition
$\mu(\As) \xi(\Bs) < \tfrac{1}{6}$ required for the convex program
(\ref{eq:sdp}) to provide exact recovery is slightly tighter than
that required for fundamental identifiability in the decomposition
problem. An important feature of our result is that it provides a
simple \emph{deterministic} condition for exact recovery. In
addition, note that the conditions only depend on the row/column
spaces of the low-rank matrix $\Bs$ and the support of the sparse
matrix $\As$, and not the singular values of $\Bs$ or the values of
the non-zero entries of $\As$. The reason for this is that the
non-zero entries of $\As$ and the singular values of $\Bs$ play no
role in the subgradient conditions with respect to the $\ell_1$ norm
and the nuclear norm.


In the sequel we discuss concrete classes of sparse and low-rank
matrices that have small $\mu$ and $\xi$ respectively.  We also show
that when the sparse and low-rank matrices $\As$ and $\Bs$ are drawn
from certain natural random ensembles, then the sufficient
conditions of Theorem~\ref{theo:main} are satisfied with high
probability; consequently, (\ref{eq:sdp}) provides exact recovery
with high probability for such matrices.


\subsection{Previous work using incoherence} The concept of
incoherence was studied in the context of recovering sparse
representations of vectors from a so-called ``overcomplete
dictionary'' \cite{Don:03}. More concretely consider a situation in
which one is given a vector formed by a sparse linear combination of
a few elements from a combined time-frequency dictionary, i.e., a
vector formed by adding a few sinusoids and a few ``spikes''; the
goal is to recover the spikes and sinusoids that compose the vector
from the infinitely many possible solutions. Based on a notion of
time-frequency incoherence, the $\ell_1$ heuristic was shown to
succeed in recovering sparse solutions \cite{Don:01}. Incoherence is
also a concept that is implicitly used in recent work under the
title of \emph{compressed sensing}, which aims to recover
``low-dimensional'' objects such as sparse vectors
\cite{Can:06,Don2:06} and low-rank matrices \cite{Rec:07,Can:08}
given incomplete observations. Our work is closer in spirit to that
in \cite{Don:03}, and can be viewed as a method to recover the
``simplest explanation'' of a matrix given an ``overcomplete
dictionary'' of sparse and low-rank matrix atoms.


\subsection{Outline} In Section~\ref{sec:app} we elaborate on the
applications mentioned previously, and discuss the implications of
our results for each of these applications. Section~\ref{sec:inc}
formally describes conditions for fundamental identifiability in the
decomposition problem based on the quantities $\xi$ and $\mu$
defined in (\ref{eq:xim}) and (\ref{eq:mum}). We also provide a
proof of the rank-sparsity uncertainty principle of
Theorem~\ref{theo:unp}. We prove Theorem~\ref{theo:main} in
Section~\ref{sec:main}, and also provide concrete classes of sparse
and low-rank matrices that satisfy the sufficient conditions of
Theorem~\ref{theo:main}. Section~\ref{sec:sim} describes the results
of simulations of our approach applied to synthetic matrix
decomposition problems. We conclude with a discussion in
Section~\ref{sec:conc}. The Appendix provides additional details and
proofs.

\section{Applications}
\label{sec:app}

In this section we describe several applications that involve
decomposing a matrix into sparse and low-rank components.

\subsection{Graphical modeling with latent variables}
We begin with a problem in statistical model selection. In many
applications large covariance matrices are approximated as low-rank
matrices based on the assumption that a small number of
\emph{latent} factors explain most of the observed statistics (e.g.,
principal component analysis). Another well-studied class of models
are those described by graphical models \cite{Lau:96} in which the
\emph{inverse} of the covariance matrix (also called the precision
or concentration or information matrix) is assumed to be
\emph{sparse} (typically this sparsity is with respect to some
graph). We describe a model selection problem involving graphical
models \emph{with} latent variables. Let the covariance matrix of a
collection of jointly Gaussian variables be denoted by $\Sigma_{(o ~
h)}$, where $o$ represents observed variables and $h$ represents
unobserved, hidden variables. The marginal statistics corresponding
to the observed variables $o$ are given by the marginal covariance
matrix $\Sigma_o$, which is simply a submatrix of the full
covariance matrix $\Sigma_{(o ~ h)}$. Suppose, however, that we
parameterize our model by the information matrix given by $K_{(o ~
h)} = \Sigma_{(o ~ h)}^{-1}$ (such a parameterization reveals the
connection to graphical models). In such a parameterization, the
\emph{marginal information matrix} corresponding to the inverse
$\Sigma_o^{-1}$ is given by the Schur complement with respect to the
block $K_h$:
\begin{equation}
\hat{K}_{o} = \Sigma_o^{-1} = K_o - K_{o,h} K_h^{-1} K_{h,o}.
\label{eq:schur}
\end{equation}
Thus if we only observe the variables $o$, we only have access to
$\Sigma_o$ (or $\hat{K}_o$). A simple explanation of the statistical
structure underlying these variables involves recognizing the
presence of the latent, unobserved variables $h$. However
(\ref{eq:schur}) has the interesting structure that $K_o$ is often
sparse due to graphical structure amongst the observed variables
$o$, while $K_{o,h} K_h^{-1} K_{h,o}$ has low-rank if the number of
latent, unobserved variables $h$ is small relative to the number of
observed variables $o$ (the rank is equal to the number of latent
variables $h$). Therefore, decomposing $\hat{K}_o$ into these sparse
and low-rank components reveals the graphical structure in the
observed variables as well as the effect due to (and the
\emph{number} of) the unobserved latent variables. We discuss this
application in more detail in a separate report \cite{Cha2:09}.


\subsection{Matrix rigidity}
The \emph{rigidity} of a matrix $M$, denoted by $R_M(k)$, is the
smallest number of entries that need to be changed in order to
reduce the rank of $M$ below $k$. Obtaining bounds on rigidity has a
number of implications in complexity theory \cite{Lok:95}, such as
the trade-offs between size and depth in arithmetic circuits.
However, computing the rigidity of a matrix is in general an NP-hard
problem \cite{Cod:00}. For any $M \in \R^{n \times n}$ one can check
that $R_M(k) \leq (n - k)^2$ (this follows directly from a Schur
complement argument). Generically every $M \in \R^{n \times n}$ is
very rigid, i.e., $R_M(k) = (n - k)^2$ \cite{Val:77}, although
special classes of matrices may be less rigid. We show that the SDP
(\ref{eq:sdp}) can be used to compute rigidity for certain matrices
with sufficiently small rigidity (see Section~\ref{subsec:rand} for
more details). Indeed, this convex program (\ref{eq:sdp}) also
provides a certificate of the sparse and low-rank components that
form such low-rigidity matrices; that is, the SDP (\ref{eq:sdp}) not
only enables us to compute the rigidity for certain matrices but
additionally provides the changes required in order to realize a
matrix of lower rank.

\subsection{Composite system identification}
A decomposition problem can also be posed in the system
identification setting. Linear time-invariant (LTI) systems can be
represented by Hankel matrices, where the matrix represents the
input-output relationship of the system \cite{Son:98}. Thus, a
sparse Hankel matrix corresponds to an LTI system with a sparse
impulse response. A low-rank Hankel matrix corresponds to a system
with small model order, and provides a minimal realization for a
system \cite{Faz:03}. Given an LTI system $H$ as follows
\begin{equation*}
H = H_{s} + H_{lr},
\end{equation*}
where $H_s$ is sparse and $H_{lr}$ is low-rank, obtaining a simple
description of $H$ requires decomposing it into its simpler sparse
and low-rank components. One can obtain these components by solving
our rank-sparsity decomposition problem. Note that in practice one
can impose in (\ref{eq:sdp}) the additional constraint that the
sparse and low-rank matrices have Hankel structure.

\subsection{Partially coherent decomposition in optical systems} We
outline an optics application that is described in greater detail in
\cite{Faz:98}. Optical imaging systems are commonly modeled using
the Hopkins integral \cite{Goo:04}, which gives the output intensity
at a point as a function of the input transmission via a quadratic
form. In many applications the operator in this quadratic form can
be well-approximated by a (finite) positive semi-definite matrix.
Optical systems described by a low-pass filter are called
\emph{coherent} imaging systems, and the corresponding system
matrices have \emph{small rank}. For systems that are not perfectly
coherent various methods have been proposed to find an \emph{optimal
coherent decomposition} \cite{Pat:94}, and these essentially
identify the best approximation of the system matrix by a matrix of
lower rank. At the other end are \emph{incoherent} optical systems
that allow some high frequencies, and are characterized by system
matrices that are \emph{diagonal}. As most real-world imaging
systems are some combination of coherent and incoherent, it was
suggested in \cite{Faz:98} that optical systems are better described
by a sum of coherent and incoherent systems rather than by the best
coherent (i.e., low-rank) approximation as in \cite{Pat:94}. Thus,
decomposing an imaging system into coherent and incoherent
components involves splitting the optical system matrix into
low-rank and diagonal components. Identifying these simpler
components has important applications in tasks such as optical
microlithography \cite{Pat:94,Goo:04}.


\section{Rank-Sparsity Incoherence} \label{sec:inc} Throughout this
paper, we restrict ourselves to square $n \times n$ matrices to avoid
cluttered notation. All our analysis extends to rectangular $n_1
\times n_2$ matrices, if we simply replace $n$ by $\max(n_1,n_2)$.


\subsection{Identifiability issues}
As described in the introduction, the matrix decomposition problem can
be fundamentally ill-posed. We describe two situations in which
identifiability issues arise. These examples suggest the kinds of
additional conditions that are required in order to ensure that there
exists a unique decomposition into sparse and low-rank matrices.

First, let $\As$ be any sparse matrix and let $\Bs = e_i e_j^T$, where
$e_i$ represents the $i$-th standard basis vector. In this case, the
low-rank matrix $\Bs$ is also very sparse, and a valid
sparse-plus-low-rank decomposition might be $\hat{A} = \As + e_i
e_j^T$ and $\hat{B} = 0$. Thus, we need conditions that ensure that
the low-rank matrix is not too sparse. One way to accomplish this is
to require that the quantity $\xi(\Bs)$ be small. As will be discussed
in Section~\ref{subsec:specialslr}), if the row and column spaces of
$\Bs$ are ``incoherent'' with respect to the standard basis, i.e., the
row/column spaces are not aligned closely with any of the coordinate
axes, then $\xi(\Bs)$ is small.

Next, consider the scenario in which $\Bs$ is any low-rank matrix
and $\As = -v e_1^T$ with $v$ being the first column of $\Bs$. Thus,
$C = \As + \Bs$ has zeros in the first column, $\rk(C) = \rk(\Bs)$,
and $C$ has the same column space as $\Bs$. Therefore, a reasonable
sparse-plus-low-rank decomposition in this case might be $\hat{B} =
\Bs + \As$ and $\hat{A} = 0$. Here $\rk(\hat{B}) = \rk(\Bs)$.
Requiring that a sparse matrix $\As$ have small $\mu(\As)$ avoids
such identifiability issues. Indeed we show in
Section~\ref{subsec:specialslr} that sparse matrices with ``bounded
degree'' (i.e., few non-zero entries per row/column) have small
$\mu$.

\subsection{Tangent-space identifiability} \label{subsec:tsiden}
We begin by describing the sets of sparse and low-rank matrices.
These sets can be considered either as differentiable manifolds
(away from their singularities) or as algebraic varieties; we
emphasize the latter viewpoint here. Recall that an algebraic
variety is defined as the zero set of a system of polynomial
equations \cite{Har:95}. The variety of rank-constrained matrices is
defined as:
\begin{equation}
\mathcal{P}(k) \triangleq \{M \in \R^{n \times n} ~ | ~ \rk(M) \leq
k\}. \label{eq:rkvar}
\end{equation}
This is an algebraic variety since it can be defined through the
vanishing of all $(k+1) \times (k+1)$ minors of the matrix $M$. The
dimension of this variety is $k(2n-k)$, and it is non-singular
everywhere except at those matrices with rank less than or equal to
$k-1$. For any matrix $M \in \R^{n \times n}$, the tangent space
$T(M)$ with respect to $\mathcal{P}(\rk(M))$ at $M$ is the span of
all matrices with either the same row-space as $M$ or the same
column-space as $M$. Specifically, let $M = U \Sigma V^T$ be a
singular value decomposition (SVD) of $M$ with $U,V \in \R^{n \times
k}$, where $\rk(M) = k$. Then we have that
\begin{equation}
T(M) = \{U X^T + Y V^T ~ | ~ X,Y \in \R^{n \times k}\}.
\label{eq:t}
\end{equation}
If $\rk(M) = k$ the dimension of $T(M)$ is $k(2n-k)$. Note that we
always have $M \in T(M)$. In the rest of this paper we view $T(M)$
as a \emph{subspace} in $\R^{n \times n}$.


Next we consider the set of all matrices that are constrained by the
size of their support. Such sparse matrices can also be viewed as
algebraic varieties:
\begin{equation}
\mathcal{S}(m) \triangleq \{M \in \R^{n \times n} ~ | ~ |\s(M)| \leq
m\}.
\end{equation}
The dimension of this variety is $m$, and it is non-singular
everywhere except at those matrices with support size less than or
equal to $m-1$. In fact $\mathcal{S}(m)$ can be thought of as a
union of ${n^2 \choose m}$ subspaces, with each subspace being
aligned with $m$ of the $n^2$ coordinate axes. For any matrix $M \in
\R^{n \times n}$, the tangent space $\Omega(M)$ with respect to
$\mathcal{S}(|\s(M)|)$ at $M$ is given by
\begin{equation}
\Omega(M) = \{N \in \R^{n \times n} ~| ~ \s(N) \subseteq \s(M) \}.
\label{eq:om}
\end{equation}
If $|\s(M)| = m$ the dimension of $\Omega(M)$ is $m$. Note again
that we always have $M \in \Omega(M)$. As with $T(M)$, we view
$\Omega(M)$ as a \emph{subspace} in $\R^{n \times n}$. Since both
$T(M)$ and $\Omega(M)$ are subspaces of $\R^{n \times n}$, we can
compare vectors in these subspaces.


Before analyzing whether $(\As,\Bs)$ can be recovered in general
(for example, using the SDP (\ref{eq:sdp})), we ask a simpler
question. Suppose that we had prior information about the tangent
spaces $\Omega(\As)$ and $T(\Bs)$, in addition to being given $C =
\As + \Bs$. Can we then \emph{uniquely} recover $(\As,\Bs)$ from
$C$? Assuming such prior knowledge of the tangent spaces is
unrealistic in practice; however, we obtain useful insight into the
kinds of conditions required on sparse and low-rank matrices for
exact decomposition. Given this knowledge of the tangent spaces, a
necessary and sufficient condition for unique recovery is that the
tangent spaces $\Omega(\As)$ and $T(\Bs)$ intersect transversally:
\begin{equation*}
\Omega(\As) \, \cap \, T(\Bs) = \{0\}.
\end{equation*}
That is, the subspaces $\Omega(\As)$ and $T(\Bs)$ have a trivial
intersection. The sufficiency of this condition for unique
decomposition is easily seen. For the necessity part, suppose for the
sake of a contradiction that a non-zero matrix $M$ belongs to
$\Omega(\As) \cap T(\Bs)$; one can add and subtract $M$ from $\As$ and
$\Bs$ respectively while still having a valid decomposition, which
violates the uniqueness requirement. The following proposition, proved
in Appendix~\ref{sec:proofs}, provides a simple condition in terms of
the quantities $\mu(\As)$ and $\xi(\Bs)$ for the tangent spaces
$\Omega(\As)$ and $T(\Bs)$ to intersect transversally.

\begin{PROP}
\label{prop:trans} Given any two matrices $\As$ and $\Bs$, we have
that
\begin{equation*}
\mu(\As) \xi(\Bs) < 1 \quad \Rightarrow \quad \Omega(\As) \, \cap \, T(\Bs) =
\{0\},
\end{equation*}
where $\xi(\Bs)$ and $\mu(\As)$ are defined in (\ref{eq:xim}) and
(\ref{eq:mum}), and the tangent spaces $\Omega(\As)$ and $T(\Bs)$
are defined in (\ref{eq:om}) and (\ref{eq:t}).
\end{PROP}

Thus, both $\mu(\As)$ and $\xi(\Bs)$ being small implies that the
tangent spaces $\Omega(\As)$ and $T(\Bs)$ intersect transversally;
consequently, we can exactly recover $(\As,\Bs)$ given $\Omega(\As)$
and $T(\Bs)$. As we shall see, the condition required in
Theorem~\ref{theo:main} (see Section~\ref{subsec:mainres}) for exact
recovery using the convex program (\ref{eq:sdp}) will be simply a mild
tightening of the condition required above for unique decomposition
given the tangent spaces.

\subsection{Rank-sparsity uncertainty principle} \label{subsec:unp}
Another important consequence of Proposition~\ref{prop:trans} is
that we have an elementary proof of the following rank-sparsity
uncertainty principle.

\begin{THEO}
\label{theo:unp}
For any matrix $M \neq 0$, we have that
\begin{equation*}
\xi(M) \mu(M) \geq 1,
\end{equation*}
where $\xi(M)$ and $\mu(M)$ are as defined in (\ref{eq:xim}) and
(\ref{eq:mum}) respectively.
\end{THEO}

\emph{Proof}: Given any $M \neq 0$ it is clear that $M \in \Omega(M)
\cap T(M)$, i.e., $M$ is an element of both tangent spaces. However
$\mu(M) \xi(M) < 1$ would imply from Proposition~\ref{prop:trans}
that $\Omega(M) \cap T(M) = \{0\}$, which is a contradiction.
Consequently, we must have that $\mu(M) \xi(M) \geq 1$. $\square$

Hence, for \emph{any} matrix $M \neq 0$ both $\mu(M)$ and $\xi(M)$
cannot be simultaneously small. Note that
Proposition~\ref{prop:trans} is an assertion involving $\mu$ and
$\xi$ for (in general) \emph{different} matrices, while
Theorem~\ref{theo:unp} is a statement about $\mu$ and $\xi$ for the
\emph{same} matrix. Essentially the uncertainty principle asserts
that no matrix can be too sparse while having ``diffuse'' row and
column spaces. An extreme example is the matrix $e_i e_j^T$, which
has the property that $\mu(e_i e_j^T) \xi(e_i e_j^T) = 1$.

\section{Exact Decomposition Using Semidefinite Programming}
\label{sec:main}


We begin this section by studying the optimality conditions of the
convex program (\ref{eq:sdp}), after which we provide a proof of
Theorem~\ref{theo:main} with simple conditions that guarantee exact
decomposition. Next we discuss concrete classes of sparse and
low-rank matrices that satisfy the conditions of
Theorem~\ref{theo:main}, and can thus be uniquely decomposed using
(\ref{eq:sdp}).

\subsection{Optimality conditions} The orthogonal projection onto the space
$\Omega(\As)$ is denoted $P_{\Omega(\As)}$, which simply sets to
zero those entries with support not inside $\s(\As)$. The subspace
orthogonal to $\Omega(\As)$ is denoted $\Omega(\As)^c$, and it
consists of matrices with complementary support, i.e., supported on
$\s(\As)^c$. The projection onto $\Omega(\As)^c$ is denoted
$P_{\Omega(\As)^c}$.

Similarly the orthogonal projection onto the space $T(\Bs)$ is denoted
$P_{T(\Bs)}$. Letting $\Bs = U \Sigma V^T$ be the SVD of $\Bs$, we
have the following explicit relation for $P_{T(\Bs)}$:
\begin{equation}
P_{T(\Bs)}(M) = P_U M + M P_V - P_U M P_V. \label{eq:pt}
\end{equation}
Here $P_U = U U^T$ and $P_V = V V^T$. The space orthogonal to
$T(\Bs)$ is denoted $T(\Bs)^\bot$, and the corresponding projection
is denoted $P_{T(\Bs)^\bot}(M)$. The space $T(\Bs)^\bot$ consists of
matrices with row-space orthogonal to the row-space of $\Bs$
\emph{and} column-space orthogonal to the column-space of $\Bs$. We
have that
\begin{equation}
P_{T(\Bs)^\bot}(M) = (I_{n \times n} - P_U) M (I_{n \times n}- P_V),
\label{eq:ptp}
\end{equation}
where $I_{n \times n}$ is the $n \times n$ identity matrix.

Following standard notation in convex analysis \cite{Roc:96}, we
denote the \emph{subgradient} of a convex function $f$ at a point
$\hat{x}$ in its domain by $\partial f(\hat{x})$. The subgradient
$\partial f(\hat{x})$ consists of all $y$ such that
\begin{equation*}
f(x) \geq f(\hat{x}) + \langle y, x - \hat{x} \rangle, ~~~ \forall
x.
\end{equation*}
From the optimality conditions for a convex program \cite{Ber:03},
we have that $(\As,\Bs)$ is an optimum of (\ref{eq:sdp}) if and only
if there exists a dual $Q \in \mathbb{R}^{n \times n}$ such that
\begin{equation}
Q \in \gamma \partial \|\As\|_1 \mathrm{~~ and ~~} Q \in
\partial \|\Bs\|_\ast. \label{eq:opt}
\end{equation}
From the characterization of the subgradient of the $\ell_1$ norm, we
have that $Q \in \gamma \partial \|\As\|_1$ if and only if
\begin{equation}
P_{\Omega(\As)}(Q) = \gamma \, \mathrm{sign}(\As), \quad
\|P_{\Omega(\As)^c}(Q)\|_\infty \leq \gamma. \label{eq:subo}
\end{equation}
Here $\mathrm{sign}(\As_{i,j})$ equals $+1$ if $\As_{i,j} > 0$, $-1$
if $\As_{i,j} < 0$, and $0$ if $\As_{i,j} = 0$. We also have that $Q
\in
\partial \|\Bs\|_\ast$ if and only if \cite{Wat:92}
\begin{equation}
P_{T(\Bs)}(Q) = U V', \quad \|P_{T(\Bs)^\bot}(Q)\| \leq 1.
\label{eq:subt}
\end{equation}
Note that these are necessary and sufficient conditions for
$(\As,\Bs)$ to be \emph{an} optimum of (\ref{eq:sdp}). The following
proposition provides sufficient conditions for $(\As,\Bs)$ to be the
\emph{unique} optimum of (\ref{eq:sdp}), and it involves a slight
tightening of the conditions (\ref{eq:opt}), (\ref{eq:subo}), and
(\ref{eq:subt}).

\begin{figure}
\begin{center}
\epsfig{file=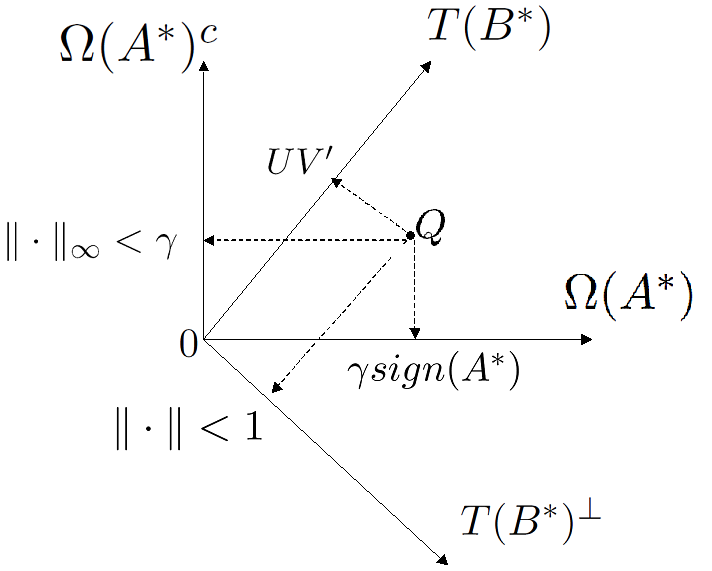,width=7cm,height=6cm}
\caption{Geometric
representation of optimality conditions: Existence of a dual $Q$. The
arrows denote orthogonal projections -- every projection must satisfy
a condition (according to Proposition~\ref{prop:sc}), which is
described next to each arrow.} 
\label{fig:fig1}
\end{center}
\end{figure}

\begin{PROP}\label{prop:sc}
Suppose that $C = \As + \Bs$. Then $(\hat{A},\hat{B}) = (\As,\Bs)$
is the \emph{unique} optimizer of (\ref{eq:sdp}) if the following
conditions are satisfied:
\begin{enumerate}
\item $\Omega(\As) \cap T(\Bs) = \{0\}$.

\item There exists a dual $Q \in \R^{n \times n}$ such that
\begin{enumerate}
\item $P_{T(\Bs)}(Q) = U V'$
\item $P_{\Omega(\As)}(Q) = \gamma \mathrm{sign}(\As)$
\item $\|P_{T(\Bs)^\bot}(Q)\| < 1$
\item $\|P_{\Omega(\As)^c}(Q)\|_{\infty} < \gamma$
\end{enumerate}
\end{enumerate}
\end{PROP}

The proof of the proposition can be found in
Appendix~\ref{sec:proofs}. Figure~\ref{fig:fig1} provides a visual
representation of these conditions. In particular, we see that the
spaces $\Omega(\As)$ and $T(\Bs)$ intersect transversely (part $(1)$
of Proposition~\ref{prop:sc}). One can also intuitively see that
guaranteeing the existence of a dual $Q$ with the requisite conditions
(part $(2)$ of Proposition~\ref{prop:sc}) is perhaps easier if the
intersection between $\Omega(\As)$ and $T(\Bs)$ is more
transverse. Note that condition $(1)$ of this proposition essentially
requires identifiability with respect to the tangent spaces, as
discussed in Section~\ref{subsec:tsiden}.

\subsection{Sufficient conditions based on $\mu(\As)$ and
$\xi(\Bs)$} \label{subsec:mainres} Next we provide simple sufficient
conditions on $\As$ and $\Bs$ that guarantee the existence of an
appropriate dual $Q$ (as required by Proposition~\ref{prop:sc}).
Given matrices $\As$ and $\Bs$ with $\mu(\As) \xi(\Bs) < 1$, we have
from Proposition~\ref{prop:trans} that $\Omega(\As) \cap T(\Bs) =
\{0\}$, i.e., condition $(1)$ of Proposition~\ref{prop:sc} is
satisfied. We prove that if a slightly stronger condition holds,
there exists a dual $Q$ that satisfies the requirements of condition
$(2)$ of Proposition~\ref{prop:sc}.

\begin{THEO}
\label{theo:main}
Given $C = \As + \Bs$ with
\begin{equation*}
\mu(\As) \xi(\Bs) < \frac{1}{6}
\end{equation*}
the \emph{unique optimum} $(\hat{A},\hat{B})$ of (\ref{eq:sdp}) is
$(\As,\Bs)$ for the following range of $\gamma$:
\begin{equation*}
\gamma \in \left(\frac{\xi(\Bs)}{1 - 4 \mu(\As) \xi(\Bs)}, \frac{1 -
3 \mu(\As) \xi(\Bs)}{\mu(\As)}\right).
\end{equation*}
Specifically $\gamma = \sqrt{\frac{3 \xi(\Bs)}{2\mu(\As)}}$ is
always inside the above range, and thus guarantees exact recovery of
$(\As,\Bs)$.

\end{THEO}


The proof of this theorem can be found in Appendix~\ref{sec:proofs}.
The main idea behind the proof is that we only consider candidates
for the dual $Q$ that lie in the direct sum $\Omega(\As) \oplus
T(\Bs)$ of the tangent spaces. Since $\mu(\As) \xi(\Bs) <
\frac{1}{6}$, we have from Proposition~\ref{prop:trans} that the
tangent spaces $\Omega(\As)$ and $T(\Bs)$ have a transverse
intersection, i.e., $\Omega(\As) \cap T(\Bs) = \{0\}$. Therefore,
there exists a \emph{unique} element $\hat{Q} \in \Omega(\As) \oplus
T(\Bs)$ that satisfies $P_{T(\Bs)}(\hat{Q}) = U V'$ and
$P_{\Omega(\As)}(\hat{Q}) = \gamma \mathrm{sign}(\As)$. The proof
proceeds by showing that if $\mu(\As) \xi(\Bs) < \frac{1}{6}$ then
the projections of this $\hat{Q}$ onto the orthogonal spaces
$\Omega(\As)^c$ and $T(\Bs)^\bot$ are small, thus satisfying
condition $(2)$ of Proposition~\ref{prop:sc}.

\paragraph{Remarks} One consequence of Theorem~\ref{theo:main} is
that if $\mu(\As) \xi(\Bs) < \frac{1}{6}$, then there exists
\emph{no} other $(A,B)$ such that $A+B = \As + \Bs$ with $\mu(A)
\xi(B) < \frac{1}{6}$. We consider this implication locally around
$(\As,\Bs)$. Recall that the quantities $\mu(\As)$ and $\xi(\Bs)$
are defined with respect to the tangent spaces $\Omega(\As)$ and
$T(\Bs)$. Suppose $\Bs$ is slightly perturbed \emph{along} the
variety of rank-constrained matrices to some $B$. This ensures that
the tangent space varies smoothly from $T(\Bs)$ to $T(B)$, and
consequently that $\xi(B) \approx \xi(\Bs)$. However, compensating
for this by changing $\As$ to $\As + (\Bs - B)$ moves $\As$ outside
the variety of sparse matrices. This is because $\Bs-B$ is \emph{not
sparse}. Thus the dimension of the tangent space $\Omega(\As + \Bs -
B)$ is much greater than that of the tangent space $\Omega(\As)$, as
a result of which $\mu(\As + \Bs - B) \gg \mu(\As)$; therefore we
have that $\xi(B) \mu(\As + \Bs - B) \gg \frac{1}{6}$. The same
reasoning holds in the opposite scenario. Consider perturbing $\As$
slightly along the variety of sparse matrices to some $A$. While
this ensures that $\mu(A) \approx \mu(\As)$, changing $\Bs$ to $\Bs
+ (\As - A)$ moves $\Bs$ outside the variety of rank-constrained
matrices. Therefore the dimension of the tangent space $T(\Bs + \As
- A)$ is greater than that of $T(\Bs)$, resulting in $\xi(\Bs + \As
- A) \gg \xi(\Bs)$; consequently we have that $\mu(A) \xi(\Bs + \As
- A) \gg \frac{1}{6}$.

\subsection{Sparse and low-rank matrices with $\mu(\As) \xi(\Bs) <
\frac{1}{6}$} \label{subsec:specialslr} We discuss concrete classes
of sparse and low-rank matrices that satisfy the sufficient
condition of Theorem~\ref{theo:main} for exact decomposition. We
begin by showing that sparse matrices with ``bounded degree'', i.e.,
bounded number of non-zeros per row/column, have small $\mu$.

\begin{PROP} \label{prop:sp}
Let $A \in \R^{n \times n}$ be any matrix with at most
$\degr_{\mathrm{max}}(A)$ non-zero entries per row/column, and with
at least $\degr_{\mathrm{min}}(A)$ non-zero entries per row/column.
With $\mu(A)$ as defined in (\ref{eq:mum}), we have that
\begin{equation*}
\degr_{\mathrm{min}}(A) \leq \mu(A) \leq \degr_{\mathrm{max}}(A).
\end{equation*}
\end{PROP}

See Appendix~\ref{sec:proofs} for the proof. Note that if $A \in
\R^{n \times n}$ has full support, i.e., $\Omega(A) = \R^{n \times
n}$, then $\mu(A) = n$. Therefore, a constraint on the number of
zeros per row/column provides a useful bound on $\mu$. We emphasize
here that simply bounding the number of non-zero entries in $A$ does
not suffice; the \emph{sparsity pattern} also plays a role in
determining the value of $\mu$.

Next we consider low-rank matrices that have small $\xi$.
Specifically, we show that matrices with row and column spaces that
are incoherent with respect to the standard basis have small $\xi$.
We measure the incoherence of a subspace $S \subseteq \R^n$ as
follows:
\begin{equation}
\beta(S) \triangleq \max_{i} \|P_{S} e_i\|_2, \label{eq:beta}
\end{equation}
where $e_i$ is the $i$'th standard basis vector, $P_S$ denotes the
projection onto the subspace $S$, and $\|\cdot\|_2$ denotes the
vector $\ell_2$ norm. This definition of incoherence also played an
important role in the results in \cite{Can:08}. A small value of
$\beta(S)$ implies that the subspace $S$ is not closely aligned with
any of the coordinate axes. In general for any $k$-dimensional
subspace $S$, we have that
\begin{equation*}
\sqrt{\frac{k}{n}} \leq \beta(S) \leq 1,
\end{equation*}
where the lower bound is achieved, for example, by a subspace that
spans any $k$ columns of an $n \times n$ orthonormal Hadamard
matrix, while the upper bound is achieved by any subspace that
contains a standard basis vector. Based on the definition of
$\beta(S)$, we define the incoherence of the row/column spaces of a
matrix $B \in \R^{n \times n}$ as
\begin{equation}
\inc(B) \triangleq \max\{\beta(\mathrm{row\mbox{-}space}(B)), ~
\beta(\mathrm{column\mbox{-}space}(B))\}. \label{eq:inc}
\end{equation}
If the SVD of $B = U \Sigma V^T$ then $\mathrm{row\mbox{-}space}(B)
= \mathrm{span}(V)$ and $\mathrm{column\mbox{-}space}(B) =
\mathrm{span}(U)$. We show in Appendix~\ref{sec:proofs} that
matrices with incoherent row/column spaces have small $\xi$; the
proof technique for the lower bound here was suggested by Ben Recht
\cite{Rec:09}.

\begin{PROP} \label{prop:lr}
Let $B \in \R^{n \times n}$ be any matrix with $\inc(B)$ defined as
in (\ref{eq:inc}), and $\xi(B)$ defined as in (\ref{eq:xim}). We
have that
\begin{equation*}
\inc(B) \leq \xi(B) \leq 2 ~ \inc(B).
\end{equation*}
\end{PROP}

If $B \in \R^{n \times n}$ is a full-rank matrix or a matrix such as
$e_1 e_1^T$, then $\xi(B) = 1$. Therefore, a bound on the
incoherence of the row/column spaces of $B$ is important in order to
bound $\xi$. Using Propositions~\ref{prop:sp} and \ref{prop:lr}
along with Theorem~\ref{theo:main} we have the following corollary,
which states that sparse bounded-degree matrices and low-rank
matrices with incoherent row/column spaces can be uniquely
decomposed.

\begin{CORL} \label{corl:deginc}
Let $C = \As + \Bs$ with $\degr_{\mathrm{max}}(\As)$ being the
maximum number of nonzero entries per row/column of $\As$ and
$\inc(\Bs)$ being the maximum incoherence of the row/column spaces
of $\Bs$ (as defined by (\ref{eq:inc})). If we have that
\begin{equation*}
\degr_{\mathrm{max}}(\As) ~ \inc(\Bs) < \frac{1}{12},
\end{equation*}
then the unique optimum of the convex program (\ref{eq:sdp}) is
$(\hat{A},\hat{B}) = (\As,\Bs)$ for a range of values of $\gamma$:
\begin{equation}
\gamma \in \left(\frac{2 ~ \inc(\Bs)}{1 - 8 ~
\degr_{\mathrm{max}}(\As) ~ \inc(\Bs)}, \frac{1 - 6 ~
\degr_{\mathrm{max}}(\As) ~
\inc(\Bs)}{\degr_{\mathrm{max}}(\As)}\right). \label{eq:grange}
\end{equation}
Specifically $\gamma = \sqrt{\frac{3 ~
\inc(\Bs)}{\degr_{\mathrm{max}}(\As)}}$ is always inside the above
range, and thus guarantees exact recovery of $(\As,\Bs)$.
\end{CORL}

We emphasize that this is a result with \emph{deterministic}
sufficient conditions on exact decomposability.

\subsection{Decomposing random sparse and low-rank matrices}
\label{subsec:rand} Next we show that sparse and low-rank matrices
drawn from certain natural random ensembles satisfy the sufficient
conditions of Corollary~\ref{corl:deginc} with high probability. We
first consider random sparse matrices with a fixed number of
non-zero entries.

\paragraph{Random sparsity model} The matrix $\As$ is such that
$\s(\As)$ is chosen uniformly at random from the collection of all
support sets of size $m$. There is no assumption made about the
values of $\As$ at locations specified by $\s(\As)$.

\begin{LEMM} \label{lemm:rsm}
Suppose that $\As \in \R^{n \times n}$ is drawn according to the
random sparsity model with $m$ non-zero entries. Let
$\degr_{\mathrm{max}}(\As)$ be the maximum number of non-zero
entries in each row/column of $\As$. We have that
\begin{equation*}
\degr_{\mathrm{max}}(\As) \leq \frac{m}{n} \log(n),
\end{equation*}
with high probability.
\end{LEMM}

The proof of this lemma follows from a standard balls and bins
argument, and can be found in several references (see for example
\cite{Bol:01}).

Next we consider low-rank matrices in which the singular vectors are
chosen uniformly at random from the set of all partial isometries.
Such a model was considered in recent work on the matrix completion
problem \cite{Can:08}, which aims to recover a low-rank matrix given
observations of a subset of entries of the matrix.

\paragraph*{Random orthogonal model \cite{Can:08}} A rank-$k$ matrix
$\Bs \in \R^{n \times n}$ with SVD $\Bs = U \Sigma V'$ is
constructed as follows: The singular vectors $U,V \in \R^{n \times
k}$ are drawn \emph{uniformly} at random from the collection of
rank-$k$ partial isometries in $\R^{n \times k}$. The choices of $U$
and $V$ need not be mutually independent. No restriction is placed
on the singular values.

As shown in \cite{Can:08}, low-rank matrices drawn from such a model
have incoherent row/column spaces.

\begin{LEMM}\label{lemm:rom}
Suppose that a rank-$k$ matrix $\Bs \in \R^{n \times n}$ is drawn
according to the random orthogonal model. Then we have that that
$\inc(\Bs)$ (defined by (\ref{eq:inc})) is bounded as
\begin{equation*}
\inc(\Bs) \lesssim \sqrt{\frac{\max(k,\log(n))}{n}},
\end{equation*}
with very high probability.
\end{LEMM}

Applying these two results in conjunction with
Corollary~\ref{corl:deginc}, we have that sparse and low-rank
matrices drawn from the random sparsity model and the random
orthogonal model can be uniquely decomposed with high probability.

\begin{CORL} \label{corl:rand}
Suppose that a rank-$k$ matrix $\Bs \in \R^{n \times n}$ is drawn
from the random orthogonal model, and that $\As \in \R^{n \times n}$
is drawn from the random sparsity model with $m$ non-zero entries.
Given $C = \As + \Bs$, there exists a range of values for $\gamma$
(given by (\ref{eq:grange})) so that $(\hat{A},\hat{B}) = (\As,\Bs)$
is the unique optimum of the SDP (\ref{eq:sdp}) with high
probability provided
\begin{equation*}
m \lesssim \frac{n^{1.5}}{\log n \sqrt{\max(k,\log n)}}.
\end{equation*}
\end{CORL}

Thus, for matrices $\Bs$ with rank $k$ smaller than $n$ the SDP
(\ref{eq:sdp}) yields exact recovery with high probability even when
the size of the support of $\As$ is super-linear in $n$. During
final preparation of this manuscript we learned of related
contemporaneous work \cite{Wri:09} that specifically studies the
problem of decomposing random sparse and low-rank matrices. In
addition to the assumptions of our random sparsity and random
orthogonal models, \cite{Wri:09} also requires that the non-zero
entries of $\As$ have independently chosen signs that are $\pm 1$
with equal probability, while the left and right singular vectors of
$\Bs$ are chosen independent of each other. For this particular
specialization of our more general framework, the results in
\cite{Wri:09} improve upon our bound in Corollary~\ref{corl:rand}.


\paragraph{Implications for the matrix rigidity problem}

Corollary~\ref{corl:rand} has implications for the matrix rigidity
problem discussed in Section~\ref{sec:app}. Recall that $R_M(k)$ is
the smallest number of entries of $M$ that need to be changed to
reduce the rank of $M$ below $k$ (the changes can be of arbitrary
magnitude). A generic matrix $M \in \R^{n \times n}$ has rigidity
$R_M(k) = (n-k)^2$ \cite{Val:77}. However, special structured
classes of matrices can have low rigidity. Consider a matrix $M$
formed by adding a sparse matrix drawn from the random sparsity
model with support size $\mathcal{O}(\frac{n}{\log n})$, and a
low-rank matrix drawn from the random orthogonal model with rank
$\epsilon n$ for some fixed $\epsilon > 0$. Such a matrix has
rigidity $R_M(\epsilon n) = \mathcal{O}(\frac{n}{\log n})$, and one
can recover the sparse and low-rank components that compose $M$ with
high probability by solving the SDP (\ref{eq:sdp}). To see this,
note that
\begin{equation*}
\frac{n}{\log n} \lesssim \frac{n^{1.5}}{\log n \sqrt{\max(\epsilon
n,\log n)}} = \frac{n^{1.5}}{\log n \sqrt{\epsilon n}},
\end{equation*}
which satisfies the sufficient condition of
Corollary~\ref{corl:rand} for exact recovery. Therefore, while the
rigidity of a matrix is NP-hard to compute in general \cite{Cod:00},
for such low-rigidity matrices $M$ one can compute the rigidity
$R_M(\epsilon n)$; in fact the SDP (\ref{eq:sdp}) provides a
certificate of the sparse and low-rank matrices that form the low
rigidity matrix $M$.


\section{Simulation Results} \label{sec:sim} We confirm the
theoretical predictions in this paper with some simple experimental
results. We also present a heuristic to choose the trade-off
parameter $\gamma$. All our simulations were performed using YALMIP
\cite{Yal:04} and the SDPT3 software \cite{Toh} for solving SDPs.

\begin{figure}
\begin{center}
\epsfig{file=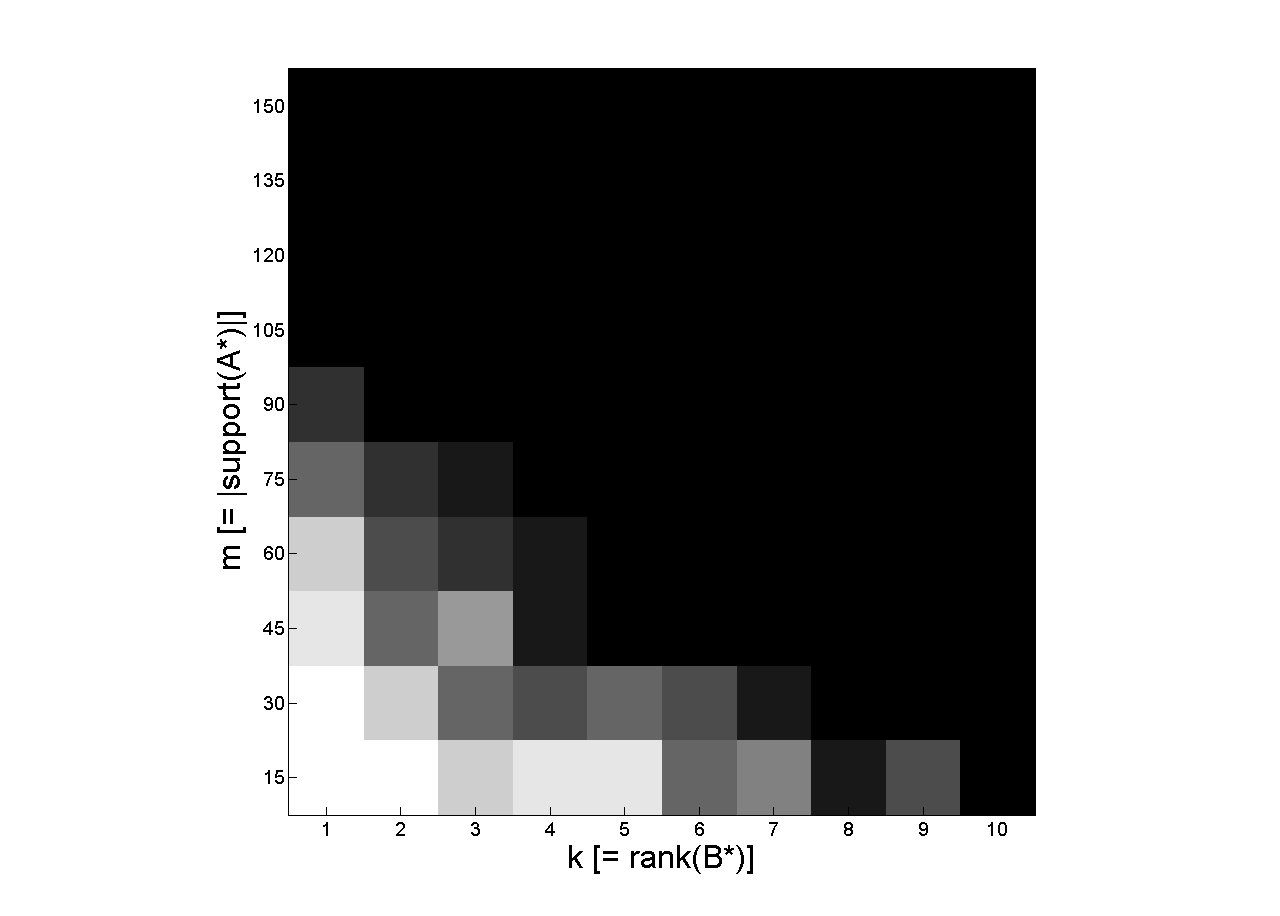,width=5.5cm} \caption{For each value of
$m,k$, we generate $25 \times 25$ random $m$-sparse $\As$ and random
rank-$k$ $\Bs$ and attempt to recover $(\As,\Bs)$ from $C = \As+\Bs$
using (\ref{eq:sdp}). For each value of $m,k$ we repeated this
procedure $10$ times. The figure shows the probability of success in
recovering $(\As,\Bs)$ using (\ref{eq:sdp}) for various values of
$m$ and $k$. White represents a probability of success of $1$, while
black represents a probability of success of $0$.}
\label{fig:simfig1}
\end{center}
\end{figure}

In the first experiment we generate random $25 \times 25$ matrices
according to the random sparsity and random orthogonal models
described in Section~\ref{subsec:rand}. To generate a random
rank-$k$ matrix $\Bs$ according to the random orthogonal model, we
generate $X, Y \in \R^{25 \times k}$ with i.i.d. Gaussian entries
and set $\Bs = X Y^T$. To generate an $m$-sparse matrix $\As$
according to the random sparsity model, we choose a support set of
size $m$ uniformly at random and the values within this support are
i.i.d. Gaussian. The goal is to recover $(\As,\Bs)$ from $C = \As +
\Bs$ using the SDP (\ref{eq:sdp}). Let $\mathrm{tol}_\gamma$ be
defined as:
\begin{equation}
\mathrm{tol}_\gamma = \frac{\|\hat{A} - \As\|_F}{\|\As\|_F} +
\frac{\|\hat{B} - \Bs\|_F}{\|\Bs\|_F}, \label{eq:tol}
\end{equation}
where $(\hat{A},\hat{B})$ is the solution of (\ref{eq:sdp}), and
$\|\cdot\|_F$ is the Frobenius norm. We declare success in
recovering $(\As,\Bs)$ if $\mathrm{tol}_\gamma < 10^{-3}$. (We
discuss the issue of choosing $\gamma$ in the next experiment.)
Figure~\ref{fig:simfig1} shows the success rate in recovering
$(\As,\Bs)$ for various values of $m$ and $k$ (averaged over $10$
experiments for each $m,k$). Thus we see that one can recover
sufficiently sparse $\As$ and sufficiently low-rank $\Bs$ from $C =
\As + \Bs$ using (\ref{eq:sdp}).

\begin{figure}
\begin{center}
\epsfig{file=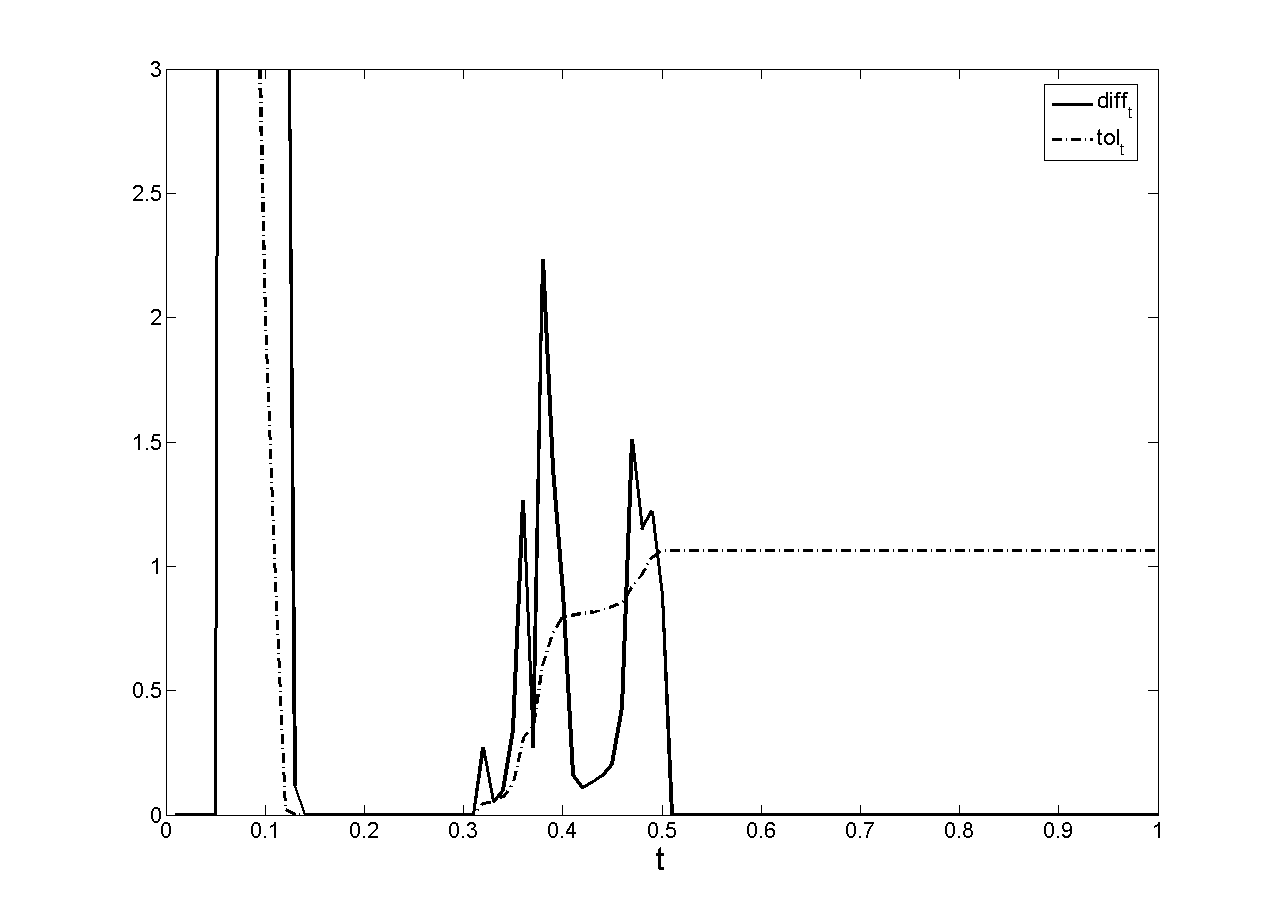,width=8cm} \caption{Comparison between
$\mathrm{tol}_t$ and $\mathrm{diff}_t$ for a randomly generated
example with $n = 25, m = 25, k = 2$.} \label{fig:simfig2}
\end{center}
\end{figure}

Next we consider the problem of choosing the trade-off parameter
$\gamma$. Based on Theorem~\ref{theo:main} we know that exact
recovery is possible for a \emph{range} of $\gamma$. Therefore, one
can simply check the stability of the solution $(\hat{A},\hat{B})$
as $\gamma$ is varied without knowing the appropriate range for
$\gamma$ in advance. To formalize this scheme we consider the
following SDP for $t \in [0,1]$, which is a slightly modified
version of (\ref{eq:sdp}):
\begin{eqnarray}
(\hat{A}_t,\hat{B}_t) = \arg \min_{A,B} &  ~~ t \|A\|_{1} +
(1-t)\|B\|_{\ast} \nonumber \\ \mbox{s.t.} & ~~ A + B = C. \label{eq:sdpt}
\end{eqnarray}
There is a one-to-one correspondence between (\ref{eq:sdp}) and
(\ref{eq:sdpt}) given by $t = \frac{\gamma}{1+\gamma}$. The benefit
in looking at (\ref{eq:sdpt}) is that the range of valid parameters
is compact, i.e., $t \in [0,1]$, as opposed to the situation in
(\ref{eq:sdp}) where $\gamma \in [0,\infty)$. We compute the
difference between solutions for some $t$ and $t-\epsilon$ as
follows:
\begin{equation}
\mathrm{diff}_t = (\|\hat{A}_{t-\epsilon}-\hat{A}_t\|_F) +
(\|\hat{B}_{t-\epsilon}-\hat{B}_t\|_F), \label{eq:diff}
\end{equation}
where $\epsilon > 0$ is some small fixed constant, say
$\epsilon=0.01$. We generate a random $\As \in \R^{25 \times 25}$
that is $25$-sparse and a random $\Bs \in \R^{25 \times 25}$ with
rank $= 2$ as described above. Given $C = \As + \Bs$, we solve
(\ref{eq:sdpt}) for various values of $t$. Figure~\ref{fig:simfig2}
shows two curves -- one is $\mathrm{tol}_t$ (which is defined
analogous to $\mathrm{tol}_\gamma$ in (\ref{eq:tol})) and the other
is $\mathrm{diff}_t$. Clearly we do not have access to
$\mathrm{tol}_t$ in practice. However, we see that $\mathrm{diff}_t$
is near-zero in exactly three regions. For sufficiently small $t$
the optimal solution to (\ref{eq:sdpt}) is $(\hat{A}_t,\hat{B}_t) =
(\As+\Bs,0)$, while for sufficiently large $t$ the optimal solution
is $(\hat{A}_t,\hat{B}_t) = (0,\As+\Bs)$. As seen in the figure,
$\mathrm{diff}_t$ stabilizes for small and large $t$. The third
``middle'' range of stability is where we typically have
$(\hat{A}_t,\hat{B}_t) = (\As,\Bs)$. Notice that outside of these
three regions $\mathrm{diff}_t$ is not close to $0$ and in fact
changes rapidly. Therefore if a reasonable guess for $t$ (or
$\gamma$) is not available, one could solve (\ref{eq:sdpt}) for a
range of $t$ and choose a solution corresponding to the ``middle''
range in which $\mathrm{diff}_t$ is stable and near zero. A related
method to check for stability is to compute the sensitivity of the
cost of the optimal solution with respect to $\gamma$, which can be
obtained from the dual solution.


\section{Discussion} \label{sec:conc} We have studied the problem of
exactly decomposing a given matrix $C = \As + \Bs$ into its sparse
and low-rank components $\As$ and $\Bs$. This problem arises in a
number of applications in model selection, system identification,
complexity theory, and optics. We characterized fundamental
identifiability in the decomposition problem based on a notion of
rank-sparsity incoherence, which relates the sparsity pattern of a
matrix and its row/column spaces via an uncertainty principle. As
the general decomposition problem is NP-hard we propose a natural
SDP relaxation (\ref{eq:sdp}) to solve the problem, and provide
sufficient conditions on sparse and low-rank matrices so that the
SDP exactly recovers such matrices. Our sufficient conditions are
deterministic in nature; they essentially require that the sparse
matrix must have support that is not too concentrated in any
row/column, while the low-rank matrix must have row/column spaces
that are not closely aligned with the coordinate axes. Our analysis
centers around studying the tangent spaces with respect to the
algebraic varieties of sparse and low-rank matrices. Indeed the
sufficient conditions for identifiability and for exact recovery
using the SDP can also be viewed as requiring that certain tangent
spaces have a transverse intersection. We also demonstrated the
implications of our results for the matrix rigidity problem.

An interesting problem for further research is the development of
special-purpose algorithms that take advantage of structure in
(\ref{eq:sdp}) to provide a more efficient solution than a
general-purpose SDP solver. Another question that arises in
applications such as model selection (due to noise or finite sample
effects) is to \emph{approximately} decompose a matrix into sparse
and low-rank components.

\section*{Acknowledgments}
The authors would like to thank Dr. Benjamin Recht and Prof. Maryam
Fazel for helpful discussions.

\appendix

\section{SDP formulation} \label{sec:sdp}
The problem (\ref{eq:sdp}) can be recast as a \emph{semidefinite
program} (SDP). We appeal to the fact that the spectral norm
$\|\cdot\|$ is the dual norm of the nuclear norm $\|\cdot\|_\ast$:
\begin{equation*}
\|M\|_\ast = \max\{\mathrm{trace}(M'Y)|~ \|Y\| \leq 1\}.
\end{equation*}
Further, the spectral norm admits a simple semidefinite
characterization \cite{Rec:07}:
\begin{equation*}
\|Y\| = \min_t ~~ t ~~~~~ \mbox{s.t.} \left(
                        \begin{array}{cc}
                          tI_{n} & Y \vspace{0.05in} \\
                          Y' & tI_{n} \\
                        \end{array}
                      \right) \succeq 0.
\end{equation*}
From duality, we can obtain the following SDP characterization of the
nuclear norm:
\begin{eqnarray*}
\|M\|_{\ast} = \min_{W_1,W_2} & ~~ \frac{1}{2}(\mathrm{trace}(W_1)+
\mathrm{trace}(W_2)) \\ \mbox{s.t.} & ~~ \left(
                        \begin{array}{cc}
                          W_1 & M \vspace{0.05in} \\
                          M' & W_2 \\
                        \end{array}
                      \right) \succeq 0.
\end{eqnarray*}

Putting these facts together, (\ref{eq:sdp}) can be rewritten as
\begin{equation}
\begin{aligned}
\min_{A,B,W_1,W_2,Z}  && \gamma \mathbf{1}_n^T Z \mathbf{1}_{n}
+ \frac{1}{2}(\mathrm{trace}(W_1)&+ \mathrm{trace}(W_2))  \\
\mbox{s.t.}&&  \left(
                        \begin{array}{cc}
                          W_1 & B \vspace{0.05in} \\
                          B' & W_2 \\
                        \end{array}
                      \right) & \succeq 0 \\
&& -Z_{i,j} \leq A_{i,j} & \leq Z_{i,j}, \quad \forall (i,j)
\\&&   A + B & = C.
\end{aligned}
\label{eq:sdp2}
\end{equation}
Here, $\mathbf{1}_{n} \in \R^{n}$ refers to the vector that has $1$ in
every entry.

\section{Proofs}
\label{sec:proofs}

\subsection*{Proof of Proposition~\ref{prop:trans}} We begin by
establishing that
\begin{equation}
\max_{N \in T(\Bs), ~ \|N\| \leq 1} \|P_{\Omega(\As)}(N)\| < 1
\quad \Rightarrow \quad
\Omega(\As) ~ \cap ~ T(\Bs) = \{0\}, \label{eq:trans1}
\end{equation}
where $P_{\Omega(\As)}(N)$ denotes the projection onto the space
$\Omega(\As)$. Assume for the sake of a contradiction that this
assertion is not true. Thus, there exists $N \neq 0$ such that $N
\in \Omega(\As) \cap T(\Bs)$. Scale $N$ appropriately such that
$\|N\| = 1$. Thus $N \in T(\Bs)$ with $\|N\| = 1$, but we also have
that $\|P_{\Omega(\As)}(N)\| = \|N\| = 1$ as $N \in \Omega(\As)$.
This leads to a contradiction.

Next, we show that
\begin{equation*}
\max_{N \in T(\Bs), ~ \|N\| \leq 1} \|P_{\Omega(\As)}(N)\| \leq
\mu(\As) \xi(\Bs),
\end{equation*}
which would allow us to conclude the proof of this proposition. We
have the following sequence of inequalities
\begin{eqnarray*}
\max_{N \in T(\Bs), ~ \|N\| \leq 1} \|P_{\Omega(\As)}(N)\| &\leq&
\max_{N \in T(\Bs), ~ \|N\| \leq 1} \mu(\As)
\|P_{\Omega(\As)}(N)\|_\infty \\ &\leq& \max_{N \in T(\Bs), ~ \|N\|
\leq 1} \mu(\As) \|N\|_\infty \\ &\leq& \mu(\As) \xi(\Bs).
\end{eqnarray*}
Here the first inequality follows from the definition (\ref{eq:mum})
of $\mu(\As)$ as $P_{\Omega(\As)}(N) \in \Omega(\As)$, the second
inequality is due to the fact that $\|P_{\Omega(\As)}(N)\|_\infty
\leq \|N\|_\infty$, and the final inequality follows from the
definition (\ref{eq:xim}) of $\xi(\Bs)$. $\square$

\subsection*{Proof of Proposition~\ref{prop:sc}} We first show that
$(\As,\Bs)$ is \emph{an} optimum of (\ref{eq:sdp}), before moving on
to showing uniqueness. Based on subgradient optimality conditions
applied at $(\As, \Bs)$, there must exist a dual $Q$ such that
\begin{equation*}
Q \in \gamma \partial\|\As\|_{1} \quad \mathrm{and} \quad Q \in
\partial\|\Bs\|_{\ast}.
\end{equation*}
The second condition in this proposition guarantees the existence of
a dual $Q$ that satisfies \emph{both} these subgradient conditions
simultaneously (see (\ref{eq:subo}) and (\ref{eq:subt})). Therefore,
we have that $(\As,\Bs)$ is \emph{an} optimum. Next we show that
under the conditions specified in the lemma, $(\As,\Bs)$ is also a
unique optimum. To avoid cluttered notation, in the rest of this
proof we let $\Omega = \Omega(\As)$, $T = T(\Bs)$, $\Omega^c(\As) =
\Omega^c$, and $T^\bot(\Bs) = T^\bot$.

Suppose that there is another feasible solution $(\As+N_A, \Bs+N_B)$
that is also a minimizer. We must have that $N_A + N_B = 0$ because
$\As+\Bs = C = (\As+N_A) + (\Bs+N_B)$. Applying the subgradient
property at $(\As,\Bs)$, we have that for \emph{any} subgradient
$(Q_A,Q_B)$ of the function $\gamma \|A\|_{1} + \|B\|_{\ast}$ (at
$(\As,\Bs)$)
\begin{equation}
\gamma \|\As+N_A\|_{1} + \|\Bs+N_B\|_{\ast} \geq \gamma \|\As\|_{1}
+ \|\Bs\|_{\ast} + \langle Q_A, N_A \rangle + \langle Q_B, N_B
\rangle. \label{eq:mainsubeq}
\end{equation}
Since $(Q_A,Q_B)$ is a subgradient of the function $\gamma \|A\|_{1}
+ \|B\|_{\ast}$ at $(\As,\Bs)$, we must have from (\ref{eq:subo})
and (\ref{eq:subt}) that
\begin{itemize}
\item $Q_A = \gamma \mathrm{sign}(\As) +
P_{\Omega^c}(Q_A)$, with $\|P_{\Omega^c}(Q_A)\|_\infty \leq \gamma$.

\item $Q_B = U V' + P_{T^\bot}(Q_B)$, with
$\|P_{T^\bot}(Q_B)\| \leq 1$.
\end{itemize}
Using these conditions we rewrite $\langle Q_A, N_A \rangle$ and
$\langle Q_B, N_B \rangle$. Based on the existence of the dual $Q$
as described in the lemma, we have that
\begin{eqnarray}
\langle Q_A, N_A \rangle &=& \langle \gamma \mathrm{sign}(\As) +
P_{\Omega^c}(Q_A), N_A \rangle \nonumber \\ &=& \langle Q -
P_{\Omega^c}(Q) + P_{\Omega^c}(Q_A), N_A \rangle \nonumber \\ &=&
\langle P_{\Omega^c}(Q_A) - P_{\Omega^c}(Q), N_A \rangle + \langle
Q, N_A \rangle, \label{eq:ipeq1}
\end{eqnarray}
where we have used the fact that $Q = \gamma \mathrm{sign}(\As) +
P_{\Omega^c}(Q)$. Similarly, we have that
\begin{eqnarray}
\langle Q_B, N_B \rangle &=& \langle U V' + P_{T^\bot}(Q_B), N_B
\rangle \nonumber \\ &=& \langle Q - P_{T^\bot}(Q) +
P_{T^\bot}(Q_B), N_B \rangle \nonumber \\ &=& \langle
P_{T^\bot}(Q_B) - P_{T^\bot}(Q), N_B \rangle + \langle Q, N_B
\rangle, \label{eq:ipeq2}
\end{eqnarray}
where we have used the fact that $Q = U V' + P_{T^\bot}(Q)$. Putting
(\ref{eq:ipeq1}) and (\ref{eq:ipeq2}) together, we have that
\begin{eqnarray}
\langle Q_A, N_A \rangle + \langle Q_B, N_B \rangle &=& \langle
P_{\Omega^c}(Q_A) - P_{\Omega^c}(Q), N_A \rangle \nonumber \\ && +
\langle P_{T^\bot}(Q_B) - P_{T^\bot}(Q), N_B \rangle \nonumber \\ &&
+ \langle Q, N_A + N_B \rangle \nonumber \\ &=& \langle
P_{\Omega^c}(Q_A) - P_{\Omega^c}(Q), N_A \rangle \nonumber \\ && +
\langle P_{T^\bot}(Q_B) - P_{T^\bot}(Q), N_B \rangle \nonumber \\
&=& \langle P_{\Omega^c}(Q_A) - P_{\Omega^c}(Q), P_{\Omega^c}(N_A)
\rangle \nonumber \\ && + \langle P_{T^\bot}(Q_B) - P_{T^\bot}(Q),
P_{T^\bot}(N_B) \rangle. \label{eq:ipeq3}
\end{eqnarray}
In the second equality, we used the fact that $N_A + N_B = 0$.

Since $(Q_A,Q_B)$ is \emph{any} subgradient, we have some freedom in
selecting $P_{\Omega^c}(Q_A)$ and $P_{T^\bot}(Q_B)$ as long as they
still satisfy the subgradient conditions
$\|P_{\Omega^c}(Q_A)\|_\infty \leq \gamma$ and $\|P_{T^\bot}(Q_B)\|
\leq 1$. We set $P_{\Omega^c}(Q_A) = \gamma
\mathrm{sign}(P_{\Omega^c}(N_A))$ so that
$\|P_{\Omega^c}(Q_A)\|_\infty = \gamma$ and $\langle
P_{\Omega^c}(Q_A), P_{\Omega^c}(N_A) \rangle = \gamma
\|P_{\Omega^c}(N_A)\|_{1}$. Letting $P_{T^\bot}(N_B) = \tilde{U}
\tilde{\Sigma} \tilde{V}^T$ be the singular value decomposition of
$P_{T^\bot}(N_B)$, we set $P_{T^\bot}(Q_B) = \tilde{U} \tilde{V}^T$
so that $\|P_{T^\bot}(Q_B)\| = 1$ and $\langle P_{T^\bot}(Q_B),
P_{T^\bot}(N_B) \rangle = \|P_{T^\bot}(N_B)\|_{\ast}$. Consequently,
we can simplify (\ref{eq:ipeq3}) as follows:
\begin{eqnarray*}
\langle Q_A, N_A \rangle + \langle Q_B, N_B \rangle &\geq& (\gamma -
\|P_{\Omega^c}(Q)\|_\infty)(\|P_{\Omega^c}(N_A)\|_{1}) \nonumber \\
&& + (1 - \|P_{T^\bot}(Q)\|)(\|P_{T^\bot}(N_B)\|_{\ast}).
\end{eqnarray*}
Since $\|P_{\Omega^c}(Q)\|_\infty < \gamma$ and $\|P_{T^\bot}(Q)\| <
1$, we have that $\langle Q_A, N_A \rangle + \langle Q_B, N_B
\rangle$ is strictly positive unless $P_{\Omega^c}(N_A) = 0$ and
$P_{T^\bot}(N_B) = 0$. (Note that if $\langle Q_A, N_A \rangle +
\langle Q_B, N_B \rangle > 0$ then $\gamma \|\As+N_A\|_{1} +
\|\Bs+N_B\|_{\ast} > \gamma \|\As\|_{1} + \|\Bs\|_{\ast}$.) However,
we have that $N_A + N_B = 0$. If $P_{\Omega^c}(N_A) =
P_{T^\bot}(N_B) = 0$, then $P_{\Omega}(N_A) + P_T(N_B) = 0$. In
other words, $P_{\Omega}(N_A) = -P_T(N_B)$. This can only be
possible if $P_{\Omega}(N_A) = P_T(N_B) = 0$ (as $\Omega \cap T =
\{0\}$), which implies that $N_A = N_B = 0$. Therefore, $\gamma
\|\As+N_A\|_{1} + \|\Bs+N_B\|_{\ast} > \gamma \|\As\|_{1} +
\|\Bs\|_{\ast}$ unless $N_A = N_B = 0$. $\square$

\subsection*{Proof of Theorem~\ref{theo:main}} As with the previous
proof, we avoid cluttered notation by letting $\Omega =
\Omega(\As)$, $T = T(\Bs)$, $\Omega^c(\As) = \Omega^c$, and
$T^\bot(\Bs) = T^\bot$. One can check that
\begin{equation}
\xi(\Bs) \mu(\As) < \frac{1}{6} \Rightarrow \frac{\xi(\Bs)}{1 - 4
\xi(\Bs) \mu(\As)} < \frac{1 - 3 \xi(\Bs) \mu(\As)}{\mu(\As)}.
\label{eq:gammarange}
\end{equation}
Thus, we show that if $\xi(\Bs) \mu(\As) < \frac{1}{6}$ then there
exists a range of $\gamma$ for which a dual $Q$ with the requisite
properties exists. Also note that plugging in $\xi(\Bs) \mu(\As) =
\tfrac{1}{6}$ in the above range gives the smaller range $(3
\xi(\Bs), \frac{1}{2 \mu(\As)})$ for $\gamma$; the geometric mean of
the extreme values gives $\gamma = \sqrt{\tfrac{3 \xi(\Bs)}{2
\mu(\As)}}$, which is always within the above range.

We aim to construct a dual $Q$ by considering candidates in the
direct sum $\Omega \oplus T$ of the tangent spaces. Since $\mu(\As)
\xi(\Bs) < \frac{1}{6}$, we can conclude from
Proposition~\ref{prop:trans} that there exists a \emph{unique}
$\hat{Q} \in \Omega \oplus T$ such that $P_\Omega(\hat{Q}) = \gamma
\mathrm{sign}(\As)$ and $P_T(\hat{Q}) = U V'$ (recall that these are
conditions that a dual must satisfy according to
Proposition~\ref{prop:sc}), as $\Omega \cap T = \{0\}$. The rest of
this proof shows that if $\mu(\As) \xi(\Bs) < \frac{1}{6}$ then the
projections of such a $\hat{Q}$ onto $T^\bot$ and onto $\Omega^c$
will be small, i.e., we show that $\|P_{\Omega^c}(\hat{Q})\|_\infty
< \gamma$ and $\|P_{T^\bot}(\hat{Q})\| < 1$.

We note here that $\hat{Q}$ can be \emph{uniquely} expressed as the
sum of an element of $T$ and an element of $\Omega$, i.e., $\hat{Q}
= Q_\Omega + Q_T$ with $Q_\Omega \in \Omega$ and $Q_T \in T$. The
uniqueness of the splitting can be concluded because $\Omega \cap T
= \{0\}$. Let $Q_\Omega = \gamma \mathrm{sign}(\As) +
\epsilon_\Omega$ and $Q_T = U V' + \epsilon_T$. We then have
\begin{equation*}
P_\Omega(\hat{Q}) = \gamma \mathrm{sign}(\As) + \epsilon_\Omega +
P_\Omega(Q_T) = \gamma \mathrm{sign}(\As) + \epsilon_\Omega +
P_\Omega(U V' + \epsilon_T).
\end{equation*}
Since $P_\Omega(\hat{Q}) = \gamma \mathrm{sign}(\As)$,
\begin{equation}
\epsilon_\Omega = - P_\Omega(U V' + \epsilon_T). \label{eq:epso}
\end{equation}
Similarly,
\begin{equation}
\epsilon_T = -P_T(\gamma \mathrm{sign}(\As) + \epsilon_\Omega).
\label{eq:epst}
\end{equation}
Next, we obtain the following bound on
$\|P_{\Omega^c}(\hat{Q})\|_\infty$:
\begin{eqnarray}
\|P_{\Omega^c}(\hat{Q})\|_\infty &=& \|P_{\Omega^c}(U V' +
\epsilon_T)\|_\infty \nonumber \\ &\leq& \|U V' +
\epsilon_T\|_\infty \nonumber \\ &\leq& \xi(\Bs) \|U V' +
\epsilon_T\| \nonumber \\ &\leq& \xi(\Bs) (1 + \|\epsilon_T\|),
\label{eq:poc}
\end{eqnarray}
where we obtain the second inequality based on the definition of
$\xi(\Bs)$ (since $U V' + \epsilon_T \in T$). Similarly, we can
obtain the following bound on $\|P_{T^\bot}(\hat{Q})\|$
\begin{eqnarray}
\|P_{T^\bot}(\hat{Q})\| &=& \|P_{T^\bot}(\gamma \mathrm{sign}(\As) +
\epsilon_\Omega)\| \nonumber \\ &\leq& \|\gamma \mathrm{sign}(\As) +
\epsilon_\Omega \| \nonumber \\ &\leq& \mu(\As)
\|\gamma \mathrm{sign}(\As) + \epsilon_\Omega\|_\infty \nonumber \\
&\leq& \mu(\As) (\gamma + \|\epsilon_\Omega\|_\infty),
\label{eq:ptperp}
\end{eqnarray}
where we obtain the second inequality based on the definition of
$\mu(\As)$ (since $\gamma \mathrm{sign}(\As) + \epsilon_\Omega \in
\Omega$). Thus, we can bound $\|P_{\Omega^c}(\hat{Q})\|_\infty$ and
$\|P_{T^\bot}(\hat{Q})\|$ by bounding $\|\epsilon_T\|$ and
$\|\epsilon_\Omega\|_\infty$ respectively (using the relations
(\ref{eq:epst}) and (\ref{eq:epso})).

By definition of $\xi(\Bs)$ and using (\ref{eq:epso}),
\begin{eqnarray}
\|\epsilon_\Omega\|_\infty &=& \|P_{\Omega}(U V' +
\epsilon_T)\|_\infty \nonumber \\ &\leq& \|U V' +
\epsilon_T\|_\infty \nonumber \\ &\leq& \xi(\Bs) \|U V' +
\epsilon_T\| \nonumber \\ &\leq& \xi(\Bs) (1 + \|\epsilon_T\|),
\label{eq:epsob1}
\end{eqnarray}
where the second inequality is obtained because $U V' + \epsilon_T
\in T$. Similarly, by definition of $\mu(\As)$ and using
(\ref{eq:epst})
\begin{eqnarray}
\|\epsilon_T\| &=& \|P_T(\gamma \mathrm{sign}(\As) +
\epsilon_\Omega)\| \nonumber \\ &\leq& 2 \|\gamma \mathrm{sign}(\As)
+ \epsilon_\Omega\| \nonumber \\ &\leq& 2 \mu(\As)
\|\gamma \mathrm{sign}(\As) + \epsilon_\Omega\|_\infty \nonumber \\
&\leq& 2 \mu(\As) (\gamma + \|\epsilon_\Omega\|_\infty),
\label{eq:epstb1}
\end{eqnarray}
where the first inequality is obtained because $\|P_T(M)\| \leq 2
\|M\|$, and the second inequality is obtained because $\gamma
\mathrm{sign}(\As) + \epsilon_\Omega \in \Omega$.


Putting
(\ref{eq:epsob1}) in (\ref{eq:epstb1}), we have that
\begin{eqnarray}
\|\epsilon_T\| &\leq& 2 \mu(\As) (\gamma + \xi(\Bs) (1 +
\|\epsilon_T\|)) \nonumber \\ \Rightarrow \|\epsilon_T\| &\leq&
\frac{2 \gamma \mu(\As) + 2 \xi(\Bs) \mu(\As)}{1 - 2 \xi(\Bs)
\mu(\As)}. \label{eq:epstb2}
\end{eqnarray}
Similarly, putting (\ref{eq:epstb1}) in (\ref{eq:epsob1}), we have
that
\begin{eqnarray}
\|\epsilon_\Omega\|_\infty &\leq& \xi(\Bs) (1 + 2 \mu(\As)(\gamma +
\|\epsilon_\Omega\|_\infty)) \nonumber \\ \Rightarrow
\|\epsilon_\Omega\|_\infty &\leq& \frac{\xi(\Bs) + 2 \gamma \xi(\Bs)
\mu(\As)}{1 - 2 \xi(\Bs) \mu(\As)}. \label{eq:epsob2}
\end{eqnarray}

We now show that $\|P_{T^\bot}(\hat{Q})\| < 1$. Combining
(\ref{eq:epsob2}) and (\ref{eq:ptperp}),
\begin{eqnarray*}
\|P_{T^\bot}(\hat{Q})\| &\leq& \mu(\As) \left(\gamma +
\frac{\xi(\Bs) + 2 \gamma \xi(\Bs) \mu(\As)}{1 - 2 \xi(\Bs)
\mu(\As)}\right)
\\ &=& \mu(\As) \left(\frac{\gamma + \xi(\Bs)}{1 - 2 \xi(\Bs)
\mu(\As)}\right) \\ &<& \mu(\As) \left(\frac{\frac{1 - 3 \xi(\Bs)
\mu(\As)}{\mu(\As)} + \xi(\Bs)}{1 - 2 \xi(\Bs) \mu(\As)}\right) \\
&=& 1,
\end{eqnarray*}
since $\gamma < \frac{1 - 3 \xi(\Bs) \mu(\As)}{\mu(\As)}$ by
assumption.

Finally, we show that $\|P_{\Omega^c}(\hat{Q})\|_\infty < \gamma$.
Combining (\ref{eq:epstb2}) and (\ref{eq:poc}),
\begin{eqnarray*}
\|P_{\Omega^c}(\hat{Q})\|_\infty &\leq& \xi(\Bs) \left(1 + \frac{2
\gamma \mu(\As) + 2 \xi(\Bs) \mu(\As)}{1 - 2 \xi(\Bs)
\mu(\As)}\right) \\ &=& \xi(\Bs) \left(\frac{1 + 2 \gamma
\mu(\As)}{1 - 2 \xi(\Bs) \mu(\As)}\right) \\ &=& \left[\xi(\Bs)
\left(\frac{1 + 2 \gamma
\mu(\As)}{1 - 2 \xi(\Bs) \mu(\As)}\right) - \gamma\right] + \gamma \\
&=& \left[\frac{\xi(\Bs) + 2 \gamma \xi(\Bs) \mu(\As) - \gamma
+ 2 \gamma \xi(\Bs) \mu(\As)}{1 - 2 \xi(\Bs) \mu(\As)}\right] + \gamma \\
&=& \left[\frac{\xi(\Bs) - \gamma (1 - 4 \xi(\Bs) \mu(\As))}{1 - 2
\xi(\Bs) \mu(\As)}\right] + \gamma \\ &<& \left[\frac{\xi(\Bs) -
\xi(\Bs)}{1 - 2 \xi(\Bs) \mu(\As)}\right] + \gamma \\ &=& \gamma.
\end{eqnarray*}
Here, we used the fact that $\frac{\xi(\Bs)}{1 - 4 \xi(\Bs)
\mu(\As)} < \gamma$ in the second inequality. $\square$


\subsection*{Proof of Proposition~\ref{prop:sp}} Based on the
Perron-Frobenius theorem \cite{Hor:90}, one can conclude that $\|P\|
\geq \|Q\|$ if $P_{i,j} \geq |Q_{i,j}|, ~ \forall ~ i,j$. Thus, we
need only consider the matrix that has $1$ in every location in the
support set $\Omega(A)$ and $0$ everywhere else. Based on the
definition of the spectral norm, we can re-write $\mu(A)$ as
follows:
\begin{equation}
\mu(A) = \max_{\|x\|_2 = 1, \|y\|_2 = 1} \sum_{(i,j) \in \Omega(A)}
x_i y_j. \label{eq:muref}
\end{equation}
Without loss of generality we restrict our attention to optima that
are achieved by element-wise non-negative vectors $x,y$.

\paragraph{Upper bound} Since the reformulation of $\mu(A)$ above involves
the maximization of a continuous function over a compact set, the
maximum is achieved at some point in the constraint set. Therefore,
we have that any optimal $(\hat{x},\hat{y})$ must satisfy the
following necessary optimality conditions: There exist Lagrange
multipliers $\lambda_1,\lambda_2$ such that
\begin{eqnarray*}
\nabla_x \left[\sum_{(i,j) \in \Omega(A)} x_i
y_j\right]_{(\hat{x},\hat{y})} &=& 2 \lambda_1 \hat{x} \\ \nabla_y
\left[\sum_{(i,j) \in \Omega(A)} x_i y_j\right]_{(\hat{x},\hat{y})}
&=& 2 \lambda_2 \hat{y}
\end{eqnarray*}
This reduces to the following system of equations:
\begin{eqnarray}
\sum_{(i,j) \in \Omega(A)} \hat{y}_j &=& 2 \lambda_1 \hat{x}_i, ~~
\forall i \label{eq:fori} \\ \sum_{(i,j) \in \Omega(A)} \hat{x}_i
&=& 2 \lambda_2 \hat{y}_j, ~~ \forall j \label{eq:forj}.
\end{eqnarray}
Multiplying the first system of equations (\ref{eq:fori})
element-wise by $\hat{x}$ and then summing, we have that
\begin{eqnarray*}
\sum_i \hat{x}_i \sum_{j:(i,j) \in \Omega(A)} \hat{y}_j &=& \sum_i
\hat{x}_i \times 2 \lambda_1 \hat{x}_i \\ \Rightarrow \sum_{(i,j)
\in \Omega(A)} \hat{x}_i \hat{y}_j &=& 2 \lambda_1.
\end{eqnarray*}
Similarly, we have that $\sum_{(i,j) \in \Omega(A)} \hat{x}_i
\hat{y}_j = 2 \lambda_2$, which implies that the Lagrange
multipliers are equal to each other and to one-half of the optimal
value attained
\begin{equation*}
2 \lambda_1 = 2 \lambda_2 = \sum_{(i,j) \in \Omega(A)} \hat{x}_i
\hat{y}_j \triangleq 2 \lambda.
\end{equation*}
We recall here that the optimal points $\hat{x},\hat{y}$ are
element-wise non-negative. Let $\sigma$ denote the element-wise sum
of the optimal points $\hat{x},\hat{y}$:
\begin{equation*}
\sigma = \sum_i \hat{x}_i + \sum_j \hat{y}_j.
\end{equation*}
Summing over all $i$ in (\ref{eq:fori}) and all $j$ in
(\ref{eq:forj}), we have that
\begin{eqnarray*}
\sum_i \sum_{j:(i,j) \in \Omega(A)} \hat{y}_j + \sum_j \sum_{i:(i,j)
\in \Omega(A)} \hat{x}_i &=& 2 \lambda \times \sigma \\ \Rightarrow
\sum_{(i,j) \in \Omega(A)} \hat{y}_j + \sum_{(i,j) \in \Omega(A)}
\hat{x}_i &=& 2 \lambda \times \sigma \\ \Rightarrow \sum_{j}
\degr_{\mathrm{max}}(A) \hat{y}_j + \sum_{i} \degr_{\mathrm{max}}(A)
\hat{x}_i &\geq& 2 \lambda \times \sigma \\ \Rightarrow
\degr_{\mathrm{max}}(A) \times \sigma &\geq& 2\lambda \times \sigma
\\ \Rightarrow \degr_{\mathrm{max}}(A) &\geq& 2 \lambda = \sum_{(i,j) \in
\Omega(A)} \hat{x}_i \hat{y}_j.
\end{eqnarray*}
Note that we used the fact that $\sigma \neq 0$. Thus, we have that
$\mu(A) \leq \degr_{\mathrm{max}}(A)$.

\paragraph{Lower bound} Now suppose that each row/column of $A$ has
\emph{at least} $\degr_{\mathrm{min}}(A)$ non-zero entries. Using
the reformulation (\ref{eq:muref}) of $\mu(A)$ above, we have that
\begin{equation*}
\mu(A) \geq \sum_{(i,j) \in \Omega(A)} \frac{1}{\sqrt{n}}
\frac{1}{\sqrt{n}} = \frac{|\s(A)|}{n} \geq \degr_{\mathrm{min}}(A).
\end{equation*}
Here we set $x = y = \frac{1}{\sqrt{n}}\mathbf{1},$ with
$\mathbf{1}$ representing the all-ones vector, as candidates in the
optimization problem (\ref{eq:muref}). $\square$

\subsection*{Proof of Proposition~\ref{prop:lr}} Let $B = U \Sigma
V^T$ be the SVD of $B$.

\paragraph{Upper bound} We can upper-bound $\xi(B)$ as follows
\begin{eqnarray*}
\xi(B) &=& \max_{M \in T(B), \|M\| \leq 1} \|M\|_\infty \\ &=&
\max_{M \in T(B), \|M\| \leq 1} \|P_{T(B)}(M)\|_\infty \\ &\leq&
\max_{\|M\| \leq 1} \|P_{T(B)}(M)\|_\infty \\ &\leq& \max_{M
\mathrm{~ unitary}} \|P_{T(B)}(M)\|_\infty \\ &\leq& \max_{M
\mathrm{~ unitary}} \|P_U M\|_\infty + \max_{M \mathrm{~ unitary}}
\|(I_{n \times n} - P_U) M P_V\|_\infty.
\end{eqnarray*}
For the second inequality, we have used the fact that the maximum of
a convex function over a convex set is achieved at one of the
extreme points of the constraint set. The unitary matrices are the
extreme points of the set of contractions (i.e., matrices with
spectral norm $\leq 1$). We have used $P_{T(B)}(M) = P_U M + M P_V -
P_U M P_V$ from (\ref{eq:pt}) in the last inequality, where $P_U = U
U^T$ and $P_V = V V^T$ denote the projections onto the spaces
spanned by $U$ and $V$ respectively.

We have the following simple bound for $\|P_U M\|_\infty$ with $M$
unitary:
\begin{eqnarray}
\max_{M \mathrm{~ unitary}} \|P_U M\|_\infty &=& \max_{M \mathrm{~
unitary}} \max_{i,j} e_i^T P_U M e_j \nonumber \\ &\leq& \max_{M
\mathrm{~ unitary}} \max_{i,j} \|P_U e_i\|_2 ~ \|M e_j\|_2 \nonumber \\
&=& \max_{i} \|P_U e_i\|_2 ~ \times ~ \max_{M \mathrm{~ unitary}}
\max_j \|M e_j\|_2 \nonumber \\ &=& \beta(U). \label{eq:pubound}
\end{eqnarray}
Here we used the Cauchy-Schwartz inequality in the second line, and
the definition of $\beta$ from (\ref{eq:beta}) in the last line.

Similarly, we have that
\begin{eqnarray}
\max_{M \mathrm{~ unitary}} \|(I_{n \times n} - P_U) M P_V \|_\infty
&=& \max_{M \mathrm{~ unitary}} \max_{i,j} e_i^T (I_{n \times n} -
P_U) M P_V e_j \nonumber \\ &\leq& \max_{M \mathrm{~ unitary}}
\max_{i,j} \|(I_{n \times n} - P_U) e_i\|_2 ~ \|M P_V e_j\|_2
\nonumber \\ &=& \max_i \|(I_{n \times n} - P_U) e_i\|_2 ~ \times ~
\max_{M \mathrm{~ unitary}} \max_{j} \|M P_V e_j\|_2 \nonumber \\
&\leq& 1 \times \max_{j} \|P_V e_j\|_2 \nonumber \\ &=& \beta(V).
\label{eq:pupvbound}
\end{eqnarray}

Using the definition of $\inc(B)$ from (\ref{eq:inc}) along with
(\ref{eq:pubound}) and (\ref{eq:pupvbound}), we have that
\begin{equation*}
\xi(B) \leq \beta(U) + \beta(V) \leq 2 ~ \inc(B).
\end{equation*}

\paragraph{Lower bound} Next we prove a lower bound on $\xi(B)$.
Recall the definition of the tangent space $T(B)$ from (\ref{eq:t}).
We restrict our attention to elements of the tangent space $T(B)$ of
the form $P_U M = U U^T M$ for $M$ unitary (an analogous argument
follows for elements of the form $P_V M$ for $M$ unitary). One can
check that
\begin{equation*}
\|P_U M\| = \max_{\|x\|_2 = 1, \|y\|_2 = 1} x^T P_U M y \leq
\max_{\|x\|_2 = 1} \|P_U x\|_2 ~ \max_{\|y\|_2 = 1} \|M y\|_2 \leq
1.
\end{equation*}
Therefore,
\begin{equation*}
\xi(B) \geq \max_{M \mathrm{~ unitary}} \|P_U M\|_\infty.
\end{equation*}
Thus, we only need to show that the inequality in line $(2)$ of
(\ref{eq:pubound}) is achieved by some unitary matrix $M$ in order
to conclude that $\xi(B) \geq \beta(U)$. Define the ``most aligned''
basis vector with the subspace $U$ as follows:
\begin{equation*}
i^\ast = \arg \max_{i} \|P_U e_i\|_2.
\end{equation*}
Let $M$ be any unitary matrix with one of its columns equal to
$\frac{1}{\beta(U)} P_U e_{i^\ast}$, i.e., a normalized version of
the projection onto $U$ of the most aligned basis vector. One can
check that such a unitary matrix achieves equality in line $(2)$ of
(\ref{eq:pubound}). Consequently, we have that
\begin{equation*}
\xi(B) \geq \max_{M \mathrm{~ unitary}} \|P_U M\|_\infty = \beta(U).
\end{equation*}
By a similar argument with respect to $V$, we have the lower bound
as claimed in the proposition. $\square$

\end{document}